\documentclass[a4paper,leqno, 12pt]{amsart}
\usepackage[ansinew]{inputenc}
\usepackage{amsmath}
\usepackage{graphics}
\usepackage{amsfonts}
\swapnumbers
\theoremstyle{plain}
\newtheorem{proclaim}{Lause}[section]

\newtheorem{theorem}[proclaim]{Theorem}
\newtheorem{lemma}[proclaim]{Lemma}
\newtheorem{corollary}[proclaim]{Corollary}
\theoremstyle{definition}
\newtheorem{definition}[proclaim]{Definition}
\newtheorem{example}[proclaim]{Example}
\theoremstyle{remark}
\newtheorem{remark}[proclaim]{Remark}

\begin{document}

\title{Optimal approximation rate of certain stochastic integrals}
\author{Heikki Sepp\"al\"a}
\address{Department of Mathematics and Statistics, P.O. Box 35 (MaD), 
FIN-40014 University of Jyv\"askyl\"a, Finland} 
\email{heikki.seppala@jyu.fi}
\maketitle

\hrule
\begin{scriptsize}
\subsection*{Abstract}
Given an increasing function $H:[0,1)\rightarrow [0,\infty)$ and 
$$
A_n(H):=\inf_{\tau\in \mathcal{T}_n}\left(\sum_{i=1}^n \int_{t_{i-1}}^{t_i} (t_i-t)H^2(t)dt\right)^{\frac{1}{2}},
$$ 
where $\mathcal{T}_n:=\{ \tau=(t_i)_{i=0}^n:\ 0=t_0<t_1<\cdots<t_n=1\}$, we characterize the property
$A_n(H)\leq \frac{c}{\sqrt{n}}$, and give conditions for $A_n(H)\leq \frac{c}{\sqrt{n^\beta}}$ and 
$A_n(H)\geq \frac{1}{c\sqrt{n^\beta}}$ for $\beta\in (0,1)$, both in terms of integrability properties of $H$. These results are 
applied to the approximation of certain stochastic integrals. 
\newline
\newline
\textit{Keywords}: Non-linear approximation; Stochastic integrals; Regular sequences  
\newline
\textit{Mathematics Subject Classification:} 41A25; 60H05
\end{scriptsize}

\smallskip
\hrule
\bigskip

\section{Introduction}

In this paper we estimate the size of the error which occurs when a stochastic integral is approximated discretely.
To explain the problem in more detail, we assume a stochastic process $X=(X_t)_{t\in [0,1]}$ such that 
$$
dX_t=\sigma(X_t)dW_t\quad \text{with } X_0\equiv x_0>0,
$$ 
where $W=(W_t)_{t\in [0,1]}$ is the standard Brownian motion, $\sigma$ satisfies certain regularity 
properties, and $(\mathcal{F}_t)_{t\in[0,1]}$ is the augmentation of the filtration generated by $W$. 
It is of interest to approximate discretely a stochastic integral, which can be written as
\begin{equation}\label{f(X)}
f(X_1)=\mathbb{E}f(X_1)+\int_0^1 \lambda_u dX_u,
\end{equation}
where $f:\mathbb{R}\rightarrow \mathbb{R}$ is a polynomially bounded, Borel measurable function, and $\lambda=(\lambda_t)_{t\in[0,1]}$ is a suitable adapted process. We approximate $f(X_1)$ by
\begin{equation}\label{approx}
\mathbb{E}f(X_1)+\sum_{i=1}^n \lambda_{t_{i-1}}(X_{t_i}-X_{t_{i-1}}),
\end{equation}
where $\tau^n:=(t_i)_{i=0}^{n}$ is a deterministic time net with $0=t_0<t_1<\cdots < t_{n-1}<t_n=1$. Using this approximation instead of the original stochastic integral, we obtain an approximation error 
\begin{equation}\label{simu}
f(X_1) - \mathbb{E}f(X_1) - \sum_{i=1}^n \lambda_{t_{i-1}}(X_{t_i}-X_{t_{i-1}}). 
\end{equation}
We are interested in the minimal quadratic error under the constraint that the time net used in the approximation has $n+1$ time points. 
According to \cite[Lemma 3.2 and its proof]{Geiss1}, this error is equivalent to
\begin{equation}
\label{a_n^X}
a_n^X(Z):= \inf_{\tau^n} a_X(Z,\tau^n),
\end{equation} 
where 
$a_X(Z, \tau^n):= \inf \left(\mathbb{E} |f(X_1)- \mathbb{E}f(X_1) -\sum_{i=1}^n v_{i-1}(X_{t_i}-X_{t_{i-1}})|^2 \right)^\frac{1}{2}$
and $Z=f(X_1)$ with the infimum taken over all sequences $v=(v_i)_{i=0}^{n-1}$ of $\mathcal{F}_{t_i}$-measurable step functions $v_i:\Omega\rightarrow \mathbb{R}$.
\bigskip

The approximation problem is of interest for at least two reasons.

\begin{enumerate}
\item[(a)] In stochastic finance one would like to replace a continuously adjusted hedging portfolio in the Black-Scholes option pricing model by a discretely adjusted one, as portfolios can be adjusted in practice only finitely many times. If we consider the quadratic error which occurs in this replacement (and which we can interpret as risk in finance), then we end up with the approximation problem described above. In this case $X=(X_t)_{t\in[0,1]}$ is an appropriate positive diffusion process, $f:(0,\infty)\rightarrow [0,\infty)$ is a payoff function of a European type option, and $\tau^n$ is the net of time points where the portfolio is rebalanced. 
\item[(b)] The approximation introduced above yields to an approximation of $\int_0^t\lambda_udX_u$ by 
$\sum_{i=1}^n \lambda_{t_{i-1}}(X_{t_i\wedge t}-X_{t_{i-1}\wedge t})$. The point is that the approximation itself is a stochastic integral, but the integrand $\lambda_u$ (which is usually hard to compute) is only computed $n$-times, whereas the increments $(X_{t_{i}\wedge t}-X_{t_{i-1}\wedge t})$ can be easily simulated (for example by an Euler scheme).
\end{enumerate}

There are several previous results concerning the error caused by the discrete approximation of stochastic integrals.
Under certain conditions on $Z$ and $\sigma$, C. and S. Geiss 
showed  that if $\tau^n=(\frac{i}{n})_{i=0}^n$ is the equidistant time net with 
cardinality $n+1$, then one has that 
$$ 
a_X(Z,\tau^n)\leq \frac{c}{\sqrt{n}} 
$$ 
if and only if $Z$ belongs to the Malliavin Sobolev space $\mathbb{D}_{1,2}$ \cite[Theorems 2.3, and 2.6]{Geiss1}. 
Furthermore, they proved that there exists a constant $c>0$ such that $a_n^X(Z)\geq\frac{1}{c\sqrt{n}}$ unless there 
are constants $c_0$ and $c_1$ such that $Z=c_0+c_1X_1$ a.s. \cite[Theorem 2.5]{Geiss1} (if such constants do exist, then 
$a_n^X(Z)=0$).  It is also known by 
\cite[Theorem 2.9]{Geiss1} and \cite[Theorem 3.2]{Geiss} that there exists a constant $c>0$ such that $a_n^X(Z)\leq\frac{c}{\sqrt{n}}$, 
if $Z$ has a certain polynomial smoothness measured by Besov spaces generated by real interpolation. In this case the rate 
$\frac{1}{\sqrt{n}}$ is obtained by adapted non--equidistant time nets. 

M. Hujo showed in \cite[Theorem 3]{Hujo} that for $X$ being the Brownian motion or the geometric Brownian motion, there 
exists random variables $Z=f(X_1)\in L_2(\Omega, \mathcal{F}, \mathbb{P})$ such that 
$$
\sup_{n\in\mathbb{N}}\sqrt{n}a_n^X(Z)=\infty,
$$ 
which means that the approximation rate is not always $\frac{1}{\sqrt{n}}$ even if the underlying 
process is the standard Brownian motion. However, there are no explicit examples of such functions.\\

These results lead us to the question to characterize those $Z=f(X_1)\in L_2(\Omega, \mathcal{F}, \mathbb{P})$ with 
$$
a_n^X(Z)\leq \frac{c}{\sqrt{n}} \quad\text{for some }c=c(Z)>0.
$$
According to Theorem \ref{a=H} below, there exists a constant $c=c(\sigma)>0$ such that 
$$
\frac{1}{c}a_X(Z, \tau)\leq \left(\sum_{i=1}^n
\int_{t_{i-1}}^{t_i}(t_i-t)H_X^2Z(t) dt\right)^{\frac{1}{2}}\leq ca_X(Z, \tau),
$$
where $H_XZ(t):=\left\Vert\left(\sigma^2\frac{\partial^2 F}{\partial x^2}\right)(t,X_t)\right\Vert_{L_2}$, $F:[0,1)\times I\rightarrow \mathbb{R}$ is given by $F(t,x)=\mathbb{E}(Z|X_t=x)$, $f$ and $X$ satisfy certain conditions and $I\subset \mathbb{R}$ depends on $X$. Moreover, Lemma \ref{Hcont} implies that $H_XZ$ is increasing so that we concentrate our investigations for some time on the quantity
$$
A_n(H) := \inf_{\tau\in \mathcal{T}_n}\left(\sum_{i=1}^n \int_{t_{i-1}}^{t_i} (t_i-t)H^2(t)dt\right)^{\frac{1}{2}},
$$
where the function $H:[0,1)\rightarrow [0,\infty)$ is increasing and
$$\mathcal{T}_n:=\{ \tau=(t_i)_{i=0}^n:\ 0=t_0<t_1<\cdots<t_n=1\}.$$ 
Our first main result, Theorem \ref{main}, says that 
$$
\inf_{\tau\in \mathcal{T}_n}\left(\sum_{i=1}^n \int_{t_{i-1}}^{t_i} (t_i-t)H^2(t)dt\right)^{\frac{1}{2}}\leq \frac{c}{\sqrt{n}},
$$
if and only if the function $H$ is integrable. 
Moreover, in Theorem \ref{main3} we give sufficient conditions for
$$
\inf_{\tau\in \mathcal{T}_n}\left(\sum_{i=1}^n \int_{t_{i-1}}^{t_i} (t_i-t)H^2(t)dt\right)^{\frac{1}{2}}\leq \frac{c}{\sqrt{n^\beta}}
$$
and 
$$
\inf_{\tau\in \mathcal{T}_n}\left(\sum_{i=1}^n \int_{t_{i-1}}^{t_i} (t_i-t)H^2(t)dt\right)^{\frac{1}{2}}\geq \frac{1}{c\sqrt{n^\beta}},
$$
where $\beta\in(0,1)$, in terms of the growth rate of $H:[0,1)\rightarrow [0,\infty)$. 

These results can be applied to the setting introduced above and also to other situations, for 
example to the quadratic approximation of multi--dimensional stochastic integrals (see \cite{Hujo2}, \cite{Temam} and \cite{Zhang}).
\bigskip

The paper is organized as follows: In section 2 we introduce the main results of the paper, their proofs can be found in Section 3. In Section 4 we apply the results of Section 2 to the $1$--dimensional stochastic setting. In particular, we give an example of random variables for which the approximation rate is $\frac{1}{c\sqrt{n^\beta}}\leq a_n^X(Z)\leq \frac{c}{\sqrt{n^\beta}}$, for $\beta\in (0,1)$ in case $X$ is the standard Brownian motion or the geometric Brownian motion. In Section 5 the results of Section 2 are applied to the approximation of $d$--dimensional stochastic integrals where the underlying diffusion might have a drift.


\section{Results}

To shorten the notation in the following, we say that $A \sim_c B$ if for constant $c\geq 1$ it holds that $\frac{1}{c}A\leq B \leq cA$ and for time net $\tau\in \mathcal{T}_n$ we define
$$||\tau||_\infty:=\max_{i\in{1,..., n}}\{t_i-t_{i-1}\}.$$

\begin{definition}\label{An}
Let $H:[0,1)\rightarrow \mathbb{R}$ be a non--negative measurable function. If $\tau=(t_i)_{i=0}^n\in\mathcal{T}_n$, then we define 
\[
\begin{cases}
A(H,\tau):=\left(\sum_{i=1}^n\int_{t_{i-1}}^{t_i}(t_i-t)H^2(t)dt\right)^\frac{1}{2},\\ 
A_n(H):=\inf_{\tau\in \mathcal{T}_n}A(H, \tau).
\end{cases}
\]

\end{definition}

\begin{definition} 
We say that an increasing function $H:[0,1)\rightarrow [0,\infty)$ belongs to the set
$\mathcal{A}$ if and only if 
$$ 
\Vert H\Vert_{\mathcal{A}}:=\sup_{n\in\mathbb{N}}\sqrt{n}A_n(H)<\infty,
$$ 
and to the set $\mathcal{H}$ if and only if 
$$ 
\Vert H\Vert_{\mathcal{H}}:=\int_0^1H(t)dt<\infty. 
$$ 
\end{definition}

\begin{theorem} \label{main} 
Let $H:[0,1)\rightarrow [0,\infty)$ be an increasing function. Then 
$$ 
\sup_{n\in\mathbb{N}}\sqrt{n}A_n(H)<\infty
$$ 
if and only if $\int_0^1H(t)dt<\infty$. In particular, one has that 
$$ 
\Vert H\Vert_{\mathcal{A}}\sim_{\sqrt{2}} \Vert H\Vert_{\mathcal{H}}. 
$$ 
\end{theorem}

\begin{remark}\label{opt.nets}
The proof of Theorem \ref{main} implies that $I:=\int_0^1H(t)dt<\infty$ gives
$$
A_n(H)\leq \frac{I}{\sqrt{n}}\quad \textit{for all } n\in\mathbb{N}. 
$$
This rate can be obtained by regular sequences (see 
\cite{Ritter} and \cite{Sacks}) generated by $H$. Regular sequences generated by $H$ are time nets 
$\tau^n=(t_i^n)_{i=0}^n$ for which 
$$
\int_0^{t_i^n} H(t) dt = \frac{i}{n}\int_0^1 H(t) dt 
$$ 
for all $i\in\{0,...,n\}$.
\end{remark}

Our second main result is

\begin{theorem} \label{main3} 
Let $H:[0,1)\rightarrow [0,\infty)$ be an increasing function and $\alpha\in(\frac{1}{2},1)$. Then one has the following:
\begin{enumerate}
\item  If there exists a constant $c_1\geq 1$  such that 
$$
H(t)\leq c_1\frac{(1-\log(1-t))^{-\alpha}}{1-t}\quad \text{for all } t\in [0,1),
$$
then
$$ 
A_n(H)\leq \frac{c}{\sqrt{n^{2\alpha-1}}} \quad \textit{for all } n\in\mathbb{N},
$$ 
where $c=c(\alpha)\geq1$.
\item If there exists $s\in[0,1)$ and a constant $c_2\geq 1$  such that 
$$
H(t)\geq \frac{1}{c_2}\frac{(1-\log(1-t))^{-\alpha}}{(1-t)}\quad \text{for all } t\in [s,1),
$$
then
$$ 
A_n(H)\geq \frac{1}{c}\frac{1}{\sqrt{n^{2\alpha-1}}} \quad \textit{for all } n\in\mathbb{N},
$$ 
where $c=c(s,\alpha,c_2)\geq 1$.\\
\end{enumerate}
\end{theorem}  

\begin{remark}
It follows from the arguments in \cite[Lemma 4.14, Proposition 4.16]{Geiss5} that if $H$ is increasing and there are 
$C\in(0,\infty)$, $\alpha\in(1,\infty)$ with 
$$ 
H(t)\leq \frac{C}{[\alpha+\log(1+\frac{1}{1-t})]^\alpha(1-t)}
$$ 
for all $t\in[0,1)$, then one has that 
\[
\sup_n\sqrt{n}A_n(H)<\infty \quad \text{for all } n\in\mathbb{N}.
\] 
\end{remark}

\begin{remark}\label{mainre}
Let $H:[0,1)\rightarrow [0,\infty)$ be a measurable function and $H^*(t):=\sup_{s\in[0,t]}H(s)<\infty$ for all $t\in [0,1)$. Then
the monotonicity properties of $A_n(\cdot)$ imply the following: 
\begin{enumerate}
\item[(1)] $\sup_{n\in\mathbb{N}}\sqrt{n}A_n(H)\leq \Vert H^*\Vert_\mathcal{H}$ as a consequence of Lemma 3.1.
\item[(2)] If $H^*(t)\leq c_1\frac{(1-\log(1-t))^{-\alpha}}{1-t}\quad \text{for all } t\in [0,1)$, then
$$
A_n(H)\leq \frac{c}{\sqrt{n^{2\alpha-1}}} \quad \textit{for all } n\in\mathbb{N}.
$$
\end{enumerate}

\end{remark}


\section{Proof}

In this chapter we prove Theorems \ref{main} and \ref{main3}. To prove Theorem \ref{main} we need two lemmas 
concerning the connection between $A_n(H)$ and $\int_0^1H(t) dt$, where $H$ is a 
non--negative and increasing function.

\begin{lemma}\label{H<C} 
Let $H:[0,T)\rightarrow [0,\infty)$, $T>0$, be an increasing function such that 
$$
I=\int_0^TH(t) dt < \infty.
$$
Then for all $n\in \mathbb{N}$ there exists a sequence $\tau^n=(t_i^n)_{i=0}^n$,
$0=t_0^n<t_1^n<\cdots<t_n^n=T$ such that 
$$
\quad \int_0^{t_i^n} H(t) dt = \frac{i}{n}I
$$
for all $i\leq n$ and for this sequence it holds that 
$$
\left(\sum_{i=1}^n\int_{t_{i-1}}^{t_i}(t_i-t)H^2(t)dt\right)^{\frac{1}{2}} \leq \frac{I}{\sqrt{n}}.
$$
\end{lemma}

\begin{proof}
The existence of the sequence $(t_i^n)_{i=0}^n$ for which 
$$
\quad \int_0^{t_i^n} H(t) dt = \frac{i}{n}I
$$
follows from the continuity of the integral. Now we have 
\[
\begin{split}
\sum_{i=1}^n \int_{t_{i-1}^n}^{t_i^n}(t_i^n-t)H^2(t) dt 
&= \sum_{i=1}^n \int_{t_{i-1}^n}^{t_i^n}[(t_i^n-t)H(t)]H(t) dt\\
&\leq \frac{I}{n}\sum_{i=1}^n \sup_{t\in [t_{i-1}^n, t_i^n)}(t_i^n-t)H(t).
\end{split} 
\]
Since $H$ is increasing, it is clear that 
$$
(t_i^n-t)H(t)\leq \int_{t}^{t_i^n}H(s) ds \leq \int_{t_{i-1}^n}^{t_i^n}H(s) ds
$$
for all $t\in [t_{i-1}^n, t_i^n)$. Hence 
$$
\sum_{i=1}^n\int_{t_{i-1}}^{t_i}(t_i-t)H^2(t)dt 
\leq \frac{I}{n}\sum_{i=1}^n \int_{t_{i-1}^n}^{t_i^n}H(t) dt = \frac{I^2}{n}.
$$

\end{proof}

\begin{lemma}\label{H<infty}
Let $H:[0,1)\rightarrow [0,\infty)$ be an increasing function. If for all 
$n\in \mathbb{N}$ there exists a time net $\tau^n=(t_i^n)_{i=0}^n\in\mathcal{T}_n$ such that 
$$
A(H,\tau^n) \leq \frac{c}{\sqrt{n}}
$$
for some fixed $c>0$, then $H$ is integrable and
$$
\int_0^1 H(t) dt \leq \sqrt{2}c.
$$
\end{lemma}

\begin{proof}
If $A(H, \tau^n)=0$, then $H\equiv 0$ and the claim is trivial. Assume then that $A(H, \tau^n)>0$, which implies that $H(t)>0$ for some $t\in [0,1)$. Let $a:=\inf\{t\in [0,1): H(t)>0\}$ and  $\tilde{\tau}^n=\{a\}\cup\{t_i^n\in\tau^n: t_i^n>a\}$. Since $H$ is positive on $(a,1)$, our 
assumption implies that $||\tilde{\tau}^n||_\infty\rightarrow 0$ as $n\rightarrow \infty$. Using the
Cauchy-Schwartz inequality and the assumption $A^2(H,\tau^n) \leq \frac{c^2}{n}$ we see that 
\begin{equation}\label{1} 
\begin{split}
\left[\sum_{i=1}^{n-1}H(t_{i-1}^n)(t_i^n-t_{i-1}^n)\right]^2 
&\leq n\sum_{i=1}^{n-1}H^2(t_{i-1}^n)(t_i^n-t_{i-1}^n)^2 \\
&\leq 2n\sum_{i=1}^n \int_{t_{i-1}^n}^{t_i^n}(t_i^n-t)H^2(t) dt\leq 2c^2.
\end{split}
\end{equation}
Let $b\in(a,1)$ and $0<\epsilon<\sqrt{c}$. Choose $n$ such that $b < t_{n-1}^n$ and
$$
\int_0^bH(t) dt < \sum_{i=1}^{n-1}H(t_{i-1}^n)(t_i^n-t_{i-1}^n)+\epsilon.
$$
(We can choose $n$ satisfying this, since the positivity of the function $H$ on the interval 
$(a,1)$ implies that $t_{n-1}^n\rightarrow 1$ and $||\tilde{\tau}^n||_\infty\rightarrow 0$ as 
$n\rightarrow \infty$.) Now \eqref{1} gives that 
$$
\int_0^bH(t) dt \leq \sqrt{2}c+\epsilon
$$
and since this is true for any $b\in (a,1)$ and any $\epsilon>0$, we finally have 
$$
\int_0^1H(t)dt \leq \sqrt{2}c.
$$
\end{proof}

\begin{proof}[Proof of Theorem \emph{\ref{main}}] 
Assume first that $H\in \mathcal{H}$. Then $I:=\int_0^1H(t)dt<\infty$ and Lemma \ref{H<C} implies
$$
\sqrt{n}A_n(H)\leq I\quad\text{for all }n\in \mathbb{N}
$$
and $\Vert H\Vert_{\mathcal{A}}\leq \Vert H\Vert_{\mathcal{H}}$.

Assume now that $H\in\mathcal{A}$, which means that 
$$
\sup_{n\in\mathbb{N}}\sqrt{n}A_n(H)<\infty.
$$ 
Lemma \ref{H<infty} implies
$$
\int_0^1H(t)dt\leq \sqrt{2}\sup_{n\in\mathbb{N}}\sqrt{n}A_n(H)
$$
and $\Vert H\Vert_{\mathcal{H}}\leq \sqrt{2}\Vert H\Vert_{\mathcal{A}}$.

The computations above imply
$$
\int_0^1H(t)dt<\infty\quad \text{if and only if} \quad\sup_{n\in\mathbb{N}}\sqrt{n}A_n(H)<\infty
$$
and that $\Vert H\Vert_\mathcal{H}\sim_{\sqrt{2}} \Vert H\Vert_\mathcal{A}$.
\end{proof}

\begin{lemma}\label{cor}
Let $H:[0,1)\rightarrow [0,\infty)$ be an increasing function. Then
$$
A_n(H)\leq \inf_{T\in(0,1)}\left[  
\frac{ \left(\int_0^TH(t)dt\right)^2}{n-1}+\int_{T}^1(1-t)H^2(t)dt\right]^\frac{1}{2}
$$
for all $n\geq 2$.
\end{lemma}

\begin{proof}
Let $T\in[0,1)$ and let $\tau^n=(t_i)_{i=0}^n\in\mathcal{T}_n$ be a time net such that $0=t_0<t_1<\cdots<t_{n-1}=T<t_n=1$ and
$$
\int_0^{t_i}H(t)dt=\frac{i}{n-1}\int_0^{T}H(t)dt\quad \text{for all }i=1,...,n-1.
$$
Using Lemma \ref{H<C} we get that
\[
\begin{split}
A^2(H,\tau^n)&=\sum_{i=1}^{n-1}\int_{t_{i-1}}^{t_i}(t_i-t)H^2(t)dt + \int_{t_{n-1}}^1(1-t)H^2(t)dt\\
&\leq\frac{\left(\int_0^{T}H(t)dt\right)^2}{n-1}+\int_{T}^1(1-t)H^2(t)dt.
\end{split}
\]
By definition, we have that $A_n(H)\leq A(H, \tau_n)$ and we are done.
\end{proof}

\begin{remark}
The best rate that Lemma \ref{cor} can give, is obtained by choosing $T$ such that
$$
\int_0^{T}H(t)dt = \sqrt{n-1}\left(\int_{T}^1(1-t)H^2(t)dt\right)^{1/2}.
$$
However, it is not known if Lemma \ref{cor} gives the optimal rate, i.e. we do not know 
whether the inequality
\begin{equation}\label{A>}
A_n^2(H)\geq \frac{1}{c}\inf_{T\in(0,1)}\left[\frac{ \left(\int_0^TH(t)dt\right)^2}{n-1}+\int_{T}^1(1-t)H^2(t)dt\right]
\end{equation}
holds. What we have is 
$$
A_n^2(H)=\inf_{T\in(0,1)}\left[A_{n-1}^2(H|[0,T])+\int_T^1(1-t)H^2(t)dt\right], 
$$
where $$A_{n-1}^2(H|[0,T]):=\inf_{0=t_0<\cdots<t_{n-1}=T}\sum_{i=1}^{n-1}\int_{t_{i-1}}^{t_i}(t_i-t)H^2(t)dt.$$ 
In order to obtain inequality \eqref{A>} we would need to know 
that there exists a constant $c>0$ such that
$$
A_{n-1}^2(H|[0,T])\geq \frac{1}{c}\frac{\left(\int_0^TH(t)dt\right)^2}{n-1},
$$ 
for all $n\geq 2$, but we do not know whether this is true.
\end{remark}

For the proof of Theorem \ref{main3}, we need the following lemmas.

\begin{lemma}\label{L1}
Let $\beta\in (0,1)$. Then there exists a constant $c>0$ such that
$$
\frac{(1-\log(1-t))^{-(1+\beta)}}{(1-t)^2}\sim_c \int_1^\infty z^{-\beta-2}(1-t)^{\frac{1}{z}-2}dz \quad \text{for all } t\in[0,1).
$$
\end{lemma}

\begin{proof}
Let $\psi_\beta(t)=\frac{(1-\log(1-t))^{-(1+\beta)}}{(1-t)^2}$ and 
$\varphi_\beta(t)= \int_1^\infty z^{-\beta-2}(1-t)^{\frac{1}{z}-2}dz$.
Choosing $x=-\frac{\log(1-t)}{z}$, we obtain 
\[
\begin{split}
\varphi_\beta(t)&= \int_1^\infty z^{-\beta-2}(1-t)^{\frac{1}{z}-2}dz\\
&=\int_{-\log(1-t)}^0\left(\frac{-\log(1-t)}{x}\right)^{-\beta-2}(1-t)^{-\frac{x}{\log(1-t)}-2}\ \frac{\log(1-t)dx}{x^2}\\
&=\frac{(-\log(1-t))^{-\beta-1}}{(1-t)^2}\int_0^{-\log(1-t)}x^\beta e^{-x}dx,
\end{split}
\]
since $(1-t)^{\frac{1}{\log(1-t)}}=[e^{\log(1-t)}]^\frac{1}{\log(1-t)}=e$.

The statement follows from
\[
\lim_{t\rightarrow 1}\frac{\varphi_\beta(t)}{\psi_\beta(t)}
=\int_0^\infty x^\beta e^{-x}dx \in (0, \infty).
\]
\end{proof}

\begin{lemma}\emph{\cite[Lemma 7]{Hujo}}\label{l7} 
Let $\theta\in[1,2)$ and $H_\theta:[0,1)\rightarrow [0,\infty)$, be given by 
$$
H_\theta(t)= \sqrt{(2-\theta)(1-t)^{-\theta}}\quad \text{for }t\in[0,1).
$$
Then 
$$
\inf_{(t_i)_{i=0}^n\in\mathcal{T}_n}\sum_{i=1}^n\int_{t_{i-1}}^{t_i}(t_i-t)H_\theta^2(t)dt\geq (\theta-1)^{n-1}
$$
for all $n\in\{1,2,...\}$.
\end{lemma}

\begin{lemma}\label{L2}
Let $H:[0,1)\rightarrow [0,\infty)$ be an increasing function and $\beta\in(0,1)$. If 
$$
H^2(t)\geq \int_1^\infty z^{-\beta-2}(1-t)^{\frac{1}{z}-2}dz \quad \text{for all } t\in[0,1),
$$
then 
$$
A_n(H)\geq \frac{1}{c_\beta\sqrt{n^\beta}} \quad \text{for all } n\in \mathbb{N}
$$
where $c_\beta=\sqrt{\beta(4^{\beta+2}+2^{\beta+2}+1)e}$.
\end{lemma}

\begin{proof}
Let  $g:[1, \infty) \times [0,1)\rightarrow (0, \infty)$ be given by 
$$g(z, t)=z^{-\beta-2}(1-t)^{\frac{1}{z}-2}.
$$
Then
$$
\frac{g(k, t)}{g(k+1, t)}=\left(1+\frac{1}{k}\right)^{\beta+2}(1-t)^{\frac{1}{k(k+1)}}\leq 2^{\beta+2}
$$ 
for all $k\geq 1$ and $t\in [0,1)$. We have
$$
\frac{d}{dz}g(z,t)=(-\log(1-t)-(2+\beta)z)\frac{(1-t)^{\frac{1}{z}-2}}{z^{\beta+4}}
$$
and it is easy to see that for any fixed $t\in [0,1)$ there exists $k_t\geq 2$ such that $g(z, t)$ is increasing for all $z\leq k_t-1$ and decreasing for all $z\geq k_t$. Hence
$$
\int_1^\infty g(z, t) dz\geq \sum_{k=1}^{k_t-2} g(k, t)+\sum_{k=k_t+1}^\infty g(k, t),
$$
where we treat an empty sum as zero. Since $g(k, t)\leq 2^{\beta+2}g(k+1, t)$ for all $k\geq 1$, we have
$$
g(k_t-1,t)+g(k_t,t)+g(k_t+1,t)\leq c_\beta g(k_t+1, t),
$$
with $c_\beta:=(4^{\beta+2}+2^{\beta+2}+1)$, and therefore 
\[
\begin{split}
\sum_{k=(k_t-1)}^\infty g(k, t)&=g(k_t-1,t)+g(k_t,t)+g(k_t+1,t)+\sum_{k=k_t+1}^\infty g(k+1, t)\\
&\leq c_\beta\sum_{k=k_t}^\infty g(k+1, t).
\end{split}
\]
This implies
\[
\begin{split}
\int_1^\infty g(z,t)dz&\geq \sum_{k=1}^{k_t-2} g(k, t)+\sum_{k=k_t}^\infty g(k+1, t)\\
&\geq \sum_{k=1}^{k_t-2} g(k, t)+ \frac{1}{c_\beta}\sum_{k=(k_t-1)}^\infty g(k, t)\\
&\geq\frac{1}{c_\beta}\sum_{k=1}^\infty g(k, t)
\end{split}
\]
for all $t\in [0,1)$. 

Let $a_k=2-\frac{1}{k}$ and $p_k=k^{-(1+\beta)}$. 
By assumption,
\[
\begin{split}
H^2(t)&\geq \int_1^\infty g(z,t)dz\\
         &\geq \frac{1}{c_\beta}\sum_{k=1}^\infty g(k, t) \\
         &   = \frac{1}{c_\beta}\sum_{k=1}^\infty \frac{1}{k^{\beta+1}}\frac{1}{k}(1-t)^{\frac{1}{k}-2} \\
         &   = \frac{1}{c_\beta}\sum_{k=1}^\infty p_k(2-a_k)(1-t)^{-a_k}.
\end{split}
\] 
Now
\[
 \begin{split}
A_n^2(H)&= \inf_{\tau\in\mathcal{T}_n} \sum_{i=1}^n\int_{t_{i-1}}^{t_i}(t_i-t)H^2(t) dt\\
&\geq \inf_{\tau\in\mathcal{T}_n} \sum_{i=1}^{n}\int_{t_{i-1}}^{t_i}(t_i-t)\frac{1}{c_\beta}\sum_{k=1}^\infty p_k(2-a_k)(1-t)^{-a_k} dt\\
&= \frac{1}{c_\beta}\inf_{\tau\in\mathcal{T}_n} \sum_{k=1}^{\infty}p_k \sum_{i=1}^n\int_{t_{i-1}}^{t_i}(t_i-t)(2-a_k)(1-t)^{-a_k}dt\\
&\geq \frac{1}{c_\beta}\sum_{k=1}^\infty p_k \inf_{\tau\in\mathcal{T}_n} \sum_{i=1}^n\int_{t_{i-1}}^{t_i}(t_i-t)(2-a_k)(1-t)^{-a_k}dt.
\end{split}
\]
To prove our claim it is enough to consider $n\geq 2$. We set 
$$
H_{a_k}(t)=\sqrt{(2-a_k)(1-t)^{-a_k}},
$$ 
and now Lemma \ref{l7} implies that
\[
\begin{split}
A_n^2(H)&\geq \frac{1}{c_\beta}\sum_{k=1}^{\infty}p_k(a_k-1)^{n-1}\\
           &=\frac{1}{c_\beta}\sum_{k=1}^{\infty}k^{-(1+\beta)}\left(1-\frac{1}{k}\right)^{n-1}\\
           &\geq \frac{1}{c_\beta e}\sum_{k=n}^{\infty}k^{-(1+\beta)}\\
           &\geq \frac{1}{c_\beta e \beta}n^{-\beta}\\
           &= \frac{1}{\tilde{c}_\beta n^\beta}, 
\end{split}
\]
where $\tilde{c}_\beta=e\beta c_\beta$. 
\end{proof}

\begin{lemma}\label{H>H+1}
Let $\beta\in (0,1)$ and $H:[0,1)\rightarrow [0,\infty)$ be an increasing function such that there exists a 
constant $c_1\geq 1$ for which
$$
A_n(H+1)\geq \frac{1}{c_1\sqrt{n^\beta}} \quad \text{for all }n\in \mathbb{N}.
$$
Then there exists a constant $c_2\geq 1$ such that
$$
A_n(H)\geq \frac{1}{c_2\sqrt{n^\beta}} \quad \text{for all }n\in \mathbb{N}.
$$
\end{lemma}

\begin{proof}
Assume first $n\geq\tilde{n}:=(2^{\beta+1}c_1^2)^\frac{1}{1-\beta}  
$. 
Then we have $\frac{1}{2c_1^2(2n)^\beta}\geq \frac{1}{n}$ and since
$$
A_{2n-1}^2(H+1)\leq 2[A_n^2(H)+ A_n^2(1)]\leq 2\left[A_n^2(H)+\frac{1}{2n}\right] \quad \text{for all }n\in \mathbb{N},
$$ 
we get 
$$
A_n^2(H)\geq \frac{1}{2c_1^2(2n-1)^\beta}-\frac{1}{2n}\geq \frac{1}{4c_1^2(2n)^\beta}=\frac{1}{\tilde{c}_2^2n^\beta}
$$
for all $n\geq \tilde{n}$, where $\tilde{c}_2=2^\frac{1+\beta}{2}\sqrt{2}c_1$.

If $n < \tilde{n}$, the computations above imply
$$
A_n^2(H)\geq A_{\lceil\tilde{n}\rceil}^2(H)\geq \frac{1}{\tilde{c}_2^2\lceil\tilde{n}\rceil^\beta}
\geq \frac{1}{c_2^2n^\beta},
$$
where $c_2=\tilde{c_2}\lceil\tilde{n}\rceil^{\frac{\beta}{2}}$ and $\lceil \tilde{n}\rceil:=\inf\{k\in\mathbb{Z}: \tilde{n}\leq k\}$.
\end{proof}

\begin{proof}[Proof of Theorem \emph{\ref{main3}}]\begin{enumerate}\item[] \end{enumerate}

(1) Let $T=1-e^{c_\alpha(n)}$, where $c_\alpha(n)=1-((1-\alpha)n^{1-\alpha}+1)^\frac{1}{1-\alpha}$. Then
\[
\begin{split}
\int_0^T H(t)dt &\leq c_1\int_0^T\frac{(1-\log(1-t))^{-\alpha}}{1-t}dt\\
&= \frac{c_1}{1-\alpha}[(1-\log(1-T))^{1-\alpha}-1]\\
&= c_1n^{1-\alpha}
\end{split}
\]
and
\[
\begin{split}
\int_T^1(1-t)H^2(t)dt &\leq c_1^2\int_T^1\frac{(1-\log(1-t))^{-2\alpha}}{1-t}dt \\
&= \frac{c_1^2}{2\alpha-1}(1-\log(1-T))^{1-2\alpha}\\
&= \frac{c_1^2}{2\alpha-1}((1-\alpha)n^{1-\alpha}+1)^{\frac{1-2\alpha}{1-\alpha}}\\
&\leq \frac{c_1^2(1-\alpha)^{\frac{1-2\alpha}{1-\alpha}}}{2\alpha-1}n^{1-2\alpha}
\end{split}
\]
and hence Lemma \ref{cor} says that, for $n\geq 2$,
\[
\begin{split}
A_n(H)&\leq  \left[\frac{1}{n-1}\left(\int_0^{T}H(t)dt\right)^2 + \int_{T}^1(1-t)H^2(t)dt\right]^{1/2}\\
&\leq \left[\frac{c_1^2}{n-1}n^{2-2\alpha}+c_1^2\tilde{c}_\alpha n^{1-2\alpha}\right]^{1/2}\\
&\leq c_1\frac{(2+\tilde{c}_\alpha)^\frac{1}{2}}{\sqrt{n^{2\alpha-1}}}, 
\end{split}
\]
where $\tilde{c}_\alpha=\frac{(1-\alpha)^{\frac{1-2\alpha}{1-\alpha}}}{2\alpha-1}$.\\
(2) Assume there exists a constant $c_2\geq 1$  such that 
$$
H(t)\geq \frac{(1-\log(1-t))^{-\alpha}}{c_2(1-t)}\quad \text{for all } t\in [s,1).
$$
Then there exists a constant $c_3\geq 1$ such that
$$
H(t)+1\geq \frac{(1-\log(1-t))^{-\alpha}}{c_3(1-t)}\quad \text{for all } t\in [0,1).
$$
If we write $\beta=2\alpha-1\in(0,1)$, Lemma \ref{L1} implies that there exists a constant $c_4\geq 1$ such that 
$$
(H(t)+1)^2\geq \frac{1}{c_4}\int_1^\infty z^{-\beta-2}(1-t)^{\frac{1}{z}-2}dz \quad \text{for all } t\in[0,1),
$$
and Lemma \ref{L2} implies that there exists $c_5\geq 1$ such that
$$
A_n(H+1)\geq \frac{1}{c_5\sqrt{n^\beta}} \quad \textit{for all } n\in\mathbb{N}.
$$
Finally, Lemma \ref{H>H+1} implies the existence of a constant $c\geq 1$ such that
$$ 
A_n(H)\geq  \frac{1}{c\sqrt{n^\beta}} \quad \textit{for all } n\in\mathbb{N}.
$$

\end{proof}


\section{Application: Optimal approximation rate of certain stochastic integrals}

Throughout the section, we assume a standard Brownian motion $W=(W_t)_{t\in[0,1]}$ on 
a stochastic basis $(\Omega, \mathcal{F}, \mathbb{P}, (\mathcal{F}_t)_{t\in[0,1]})$, where $(\mathcal{F}_t)_{t\in[0,1]}$
is the augmentation of the natural filtration of $W$ and $\mathcal{F}=\mathcal{F}_1$. We let the process $S=(S_t)_{t\in[0,1]}$ be 
the geometric Brownian motion, i.e. $S_t=e^{W_t-\frac{t}{2}}$ for all $t\in [0,1]$. The space of continuous, infinitely many times continuously differentiable functions with bounded derivatives is denoted by $\mathcal{C}_b^\infty(\mathbb{R})$. Moreover, we let $X=(X_t)_{t\in[0,1]}$ 
be a diffusion such that
\begin{equation}\label{X_t}
dX_t=\sigma(X_t)dW_t\quad \text{with } X_0\equiv x_0\in\mathbb{R},
\end{equation} 
where the process $X$ is obtained through $Y=(Y_t)_{t\in[0,1]}$ given as unique continuous solution of 
$$
dY_t=\hat{\sigma}(Y_t)dW_t+\hat{b}(Y_t)dt \quad \text{with } Y_0\equiv
y_0\in\mathbb{R}, 
$$
with $0<\epsilon_0\leq \hat{\sigma}\in \mathcal{C}_b^\infty(\mathbb{R})$ and 
$\hat{b}\in \mathcal{C}_b^\infty(\mathbb{R})$, in the following two ways:
\begin{enumerate}
\item[(a)] $y_0=x_0\in \mathbb{R}$, $\hat{\sigma}:=\sigma$, $\hat{b}:=0$, $X_t:=Y_t$, 
\item[(b)] $y_0=\log x_0$ with $x_0>0$, 
$$
\hat{\sigma}(y):=\frac{\sigma(e^y)}{e^y}, \quad
\hat{b}(y):=-\frac{1}{2}\hat{\sigma}(y)^2, 
\quad{and }\ X_t=e^{Y_t}.
$$ 
\end{enumerate}

Moreover, we let $\gamma$ be the Gaussian measure on $\mathbb{R}$, i.e.
$$
d\gamma(x):=\frac{1}{\sqrt{2\pi}}e^{-\frac{x^2}{2}}dx.
$$

\begin{definition} 
Let $\mathcal{C}_e$ be the linear space of Borel measurable functions 
$f:\mathbb{R}\rightarrow \mathbb{R}$ such that there exists $m>0$ for which 
$$ 
\sup_{x\in\mathbb{R}}e^{-m|x|}\mathbb{E}f^2(x+tg)<\infty 
$$ 
for all $t>0$, where $g$ is a centered standard normal random variable. Moreover, we define 
$$
\mathcal{C}:=\{Z:=f(Y_1):\Omega \rightarrow \mathbb{R}\ |\ f\in \mathcal{C}_e \text{ and } Y \text{ as above} \}.
$$
\end{definition}

The main tool for investigating the approximation problem in papers of C. Geiss, S. Geiss, and Hujo was the 
$H$-functional defined in the following way.

\begin{definition} 
Let $X$ be a stochastic process as in \eqref{X_t} and assume that $Z\in \mathcal{C}$ 
(or $Z\in L_2(\Omega, \mathcal{F}, \mathbb{P})$ if $X\in\{W,S\}$). Then we set
\begin{equation}\label{H} 
H_XZ(t):=\left\Vert\left(\sigma^2\frac{\partial^2 F}{\partial x^2}\right)(t,X_t)\right\Vert_{L_2}
\quad\text{ for all }t\in[0,1),
\end{equation}
where $F:[0,1)\times I\rightarrow \mathbb{R}$ is given by $F(t,x)=\mathbb{E}(Z|X_t=x)$, with $I=\mathbb{R}$ in the 
case of (a) and $I=(0,\infty)$ in the case of (b).
\end{definition}

\begin{lemma}\emph{\cite[Lemma 5.3]{Geiss1}, \cite[Lemma 3.9]{Geiss}} \label{Hcont}
The function $H_XZ:[0,1)\rightarrow [0,\infty)$ is continuous and increasing. 
\end{lemma}

\bigskip

In order to deduce from Theorem \ref{main} a characterization of the approximation rate 
$$
a_n^X(Z)\leq \frac{c}{\sqrt{n}},
$$
we need the following theorem.

\bigskip

\begin{theorem}\emph{\cite[Lemma 3.2]{Geiss1} \cite[Lemma 3.10]{Geiss}} \label{a=H}
Let $X$ be a stochastic process as in \eqref{X_t}, $Z\in\mathcal{C}$  
(or $Z\in L_2(\Omega, \mathcal{F}, \mathbb{P})$ if $X\in\{W,S\}$) and 
$\tau=(t_i)_{i=0}^n\in\mathcal{T}_n$. Then  
$$
a_X(Z, \tau)\sim_{c} \left(\sum_{i=1}^n
\int_{t_{i-1}}^{t_i}(t_i-t)H_X^2Z(t) 
dt\right)^{\frac{1}{2}}
$$
where  $c\geq 1$ is an absolute constant depending on $\sigma$ only. Consequently, 
$$
a_n^X(Z)\sim_{c}A_n(H_XZ).
$$
\end{theorem}

\begin{corollary} \label{maina} 
Let $X$ be as in \eqref{X_t} and $Z\in \mathcal{C}$ (or $Z\in L_2(\Omega, \mathcal{F}, \mathbb{P})$ if $X\in\{W,S\}$). Then 
$$ 
\sup_{n\in\mathbb{N}}\sqrt{n}a_n^X(Z)\sim_{c} \int_0^1\left\Vert\left(\sigma^2\frac{\partial^2 F}{\partial x^2}\right)(t,X_t)\right\Vert_{L_2}dt,
$$ 
where $F:[0,1)\times I\rightarrow \mathbb{R}$ is given by $F(t,x)=\mathbb{E}(Z|X_t=x)$, with $I=\mathbb{R}$ in the 
case of (a) and $I=(0,\infty)$ in the case of (b).
\end{corollary}  

\begin{proof}
Theorem \ref{main} together with Lemma \ref{Hcont} and Theorem \ref{a=H} gives the result immediately.
\end{proof}

\begin{remark}
Remark \ref{opt.nets} implies that if $\left\Vert\left(\sigma^2\frac{\partial^2 F}{\partial x^2}\right)(t,X_t)\right\Vert_{L_2}$ 
is integrable, then the regular sequences generated by 
$\left\Vert\left(\sigma^2\frac{\partial^2 F}{\partial x^2}\right)(t,X_t)\right\Vert_{L_2}$ give the rate $\frac{1}{\sqrt{n}}$. 
Using these sequences, denoted by $\tau_r^n$, we have that if 
$A:=\int_0^1\left\Vert\left(\sigma^2\frac{\partial^2 F}{\partial x^2}\right)(t,X_t)\right\Vert_{L_2}dt<\infty$, then  
$$
a_n^X(Z)\leq a_X(Z, \tau_r^n)\leq \frac{c_{(\ref{a=H})} A}{\sqrt{n}}\quad \textit{for all } n\in\mathbb{N}, 
$$
where $c_{(\ref{a=H})}>0$ is taken from Theorem \ref{a=H} above.

One can also optimize over random time nets instead of deterministic ones considered here. The result 
\cite[Theorem 1.1.]{Geiss3} from C. and S. Geiss implies that $\frac{1}{\sqrt{n}}$ is the best possible approximation 
rate also for the random time nets in case the underlying diffusion $X$ is the Brownian motion $W$ 
or the geometric Brownian motion $S$ and $Z$ is not equal to $c_0+c_1X_1$ a.s. for some $c_0, c_1\in\mathbb{R}$. 
This means that if $X\in\{W, S\}$, the random time nets do not improve the approximation if the deterministic 
time nets already give the rate $\frac{1}{\sqrt{n}}$. According to this, Corollary \ref{maina} 
implies that if 
$$
\int_0^1 \left\Vert\left(\sigma^2\frac{\partial^2 F}{\partial x^2}\right)(t,X_t)\right\Vert_{L_2}dt<\infty,
$$ 
then the optimal approximation rate is $\frac{1}{\sqrt{n}}$ also for the random time nets and this rate is obtained by 
using the regular sequences generated by $\left\Vert\left(\sigma^2\frac{\partial^2 F}{\partial x^2}\right)(t,X_t)\right\Vert_{L_2}$.

\end{remark}

Now we give for $\beta \in (0,1)$ an example such that 
$$
a_n^X(Z)\sim_c\frac{1}{\sqrt{n^\beta}}\quad \text{for all } n\in\mathbb{N},
$$
in case $X$ is a standard Brownian motion or the geometric Brownian motion. According to Theorem \ref{main3}, Lemma \ref{Hcont} 
and Theorem \ref{a=H} it is sufficient to find a random variable $Z=f_\alpha(W_1)$ such that 
$$
H_XZ(t)\sim_c \frac{(1-\log(1-t))^{-\alpha}}{1-t},
$$ 
where $\alpha=\frac{\beta+1}{2}$.

\begin{example}\label{mainex}
Let $\alpha\in (1/2, 1)$ and $f_\alpha=\sum_{k=0}^\infty a_kh_k\in L_2(\gamma)$, where $a=(a_k)_{k=0}^\infty$ is given by
\[
a_k=
\begin{cases}
0\quad &\text{if } k\in\{0,1,3\},\\
\frac{1}{\sqrt{2}}&\text{if } k=2,\\
\sqrt{\frac{k-2}{k(k-1)}}\log^{-\alpha}(k-2)\quad &\text{if } k\geq 4,
\end{cases}
\]
and $(h_k)_{k=0}^\infty\subset L_2(\gamma)$ is the complete orthonormal system of Hermite polynomials, 
$$
h_k(x)=\frac{(-1)^k}{\sqrt{k!}}e^{\frac{x^2}{2}}\frac{d^k}{dx^k}e^{-\frac{x^2}{2}}.
$$
Then $Z_\alpha:=f_\alpha(W_1)\in L_2(\Omega, \mathcal{F}, \mathbb{P})$ and it can be shown that
$$
H_WZ_\alpha(t)=\left(1+\sum_{k=2}^\infty k\log^{-2\alpha}(k)t^k\right)^{1/2}\sim_{c_1} \frac{(1-\log(1-t))^{-\alpha}}{1-t}
$$ 
for all $t\in [0,1)$ (according to Lemmas \ref{L3.7} and \ref{logser} below). Using Lemma \ref{L3.7} it is easy to show that 
there exists a constant $c_2>0$ such that
$$
H_WZ_\alpha(t)\sim_{c_2} H_SZ_\alpha(t) \quad \textit{for all } t\in (0,1).
$$
Theorem \ref{main3} implies there exists a constant $c_3\geq 1$ such that
\[
\frac{1}{c_3\sqrt{n^{2\alpha-1}}}\leq a_n^X(Z_\alpha)\leq \frac{c_3}{\sqrt{n^{2\alpha-1}}}
\]
for all $n\in\mathbb{N}$, where $X\in\{W, S\}$. In other words, letting $\beta\in(0,1)$ and defining $\alpha:=\frac{\beta+1}{2}$
we have
$$
a_n^X(Z_\alpha)\sim_{c_3} \frac{1}{\sqrt{n^\beta}} \quad \text{for all }n\in\mathbb{N}.
$$

\end{example}

The following lemma should be known. For completeness and convenience of the reader we include a proof.

\begin{lemma}\label{logser}
Let $\beta>1$. Then for all $t\in [0,1)$, one has that
\begin{equation}\label{logser/1}
\frac{(1-\log(1-t))^{-\beta}}{(1-t)^2}\sim_c 1+\sum_{k=2}^\infty k\log^{-\beta}(k)t^k,
\end{equation}
where the constant $c\geq 1$ depends at most on $\beta$.
\end{lemma}

\begin{proof}
Let $n\geq e^{\beta}$ be an integer, $\epsilon \in [\frac{1}{n+1}, \frac{1}{n})$, and $t=e^{-\epsilon}$. Since $k\log^{-\beta}(k)$ is increasing if $k\geq e^{\beta}$ and we assumed that $n\geq e^{\beta}$, we have 
\[
\begin{split}
1+\sum_{k=2}^\infty k\log^{-\beta}(k)t^k &\geq \sum_{k=n}^{2n} k\log^{-\beta}(k)(e^{-1/n})^k\\
                                      &\geq \sum_{k=n}^{2n} n \log^{-\beta}(n) e^{-2} \\
                                      &\geq e^{-2} n^2\log^{-\beta}(n).
\end{split}
\]

Moreover, 
\[
\begin{split}
1+\sum_{k=2}^\infty k\log^{-\beta}(k)t^k 
&\leq 1+\sum_{k=2}^n k\log^{-\beta}(k)
+\sum_{m=1}^\infty\sum_{k=mn+1}^{(m+1)n}k\log^{-\beta}(k)e^{-\frac{mn}{n+1}} \\
&\leq c_\beta\sum_{k=2}^n n\log^{-\beta}(n)
+\sum_{m=1}^\infty (m+1)n^2\log^{-\beta}(n)e^{-\frac{mn}{n+1}} \\
&\leq c_\beta n^2\log^{-\beta}(n) + n^2 \log^{-\beta}(n)\sum_{m=1}^\infty (m+1)e^{-m/2}\\ 
&\leq (c_\beta+c) n^2\log^{-\beta}(n),
\end{split} 
\] 
where $c_\beta$ depends at most on $\beta$ and $c=\sum_{m=1}^\infty (m+1)e^{-m/2}$. 
This implies, for $t=e^{-\epsilon}$ with $\epsilon\in [\frac{1}{n+1},\frac{1}{n})$, that  
$$
1+\sum_{k=2}^\infty k\log^{-\beta}(k)t^k \sim_{c_1} n^2\log^{-\beta}(n)\quad\text{for all }n\geq e^\beta, 
$$
where $c_1\geq 1$ is a constant depending at most on $\beta$. Adapting the constant $c_1>0$, we get this for $n\geq 2$. \\
Now we show that if $n\geq 4$, then
$$
\frac{(1-\log(1-t))^{-\beta}}{(1-t)^2}\sim_{c_2} n^2\log^{-\beta}(n),
$$
where $c_2\geq 2$ is a constant depending at most on $\beta$. Firstly, we have that $\log(\frac{1}{t})\sim_{c_3}\frac{1}{n}$, where $c_3=\frac{5}{4}$. Moreover
$$
\log(u^{-1})\sim_{c_4} 1-u ,
$$
for all  $u\in[e^{-1/2},1]$, where $c_4=[2(1-e^{-\frac{1}{2}})]^{-1}$. Hence 
$$
1-t \sim_{c_5} \frac{1}{n},
$$
where $c_5=\frac{5}{8}[1-e^{-\frac{1}{2}}]^{-1}$. Furthermore, 
$$
\frac{\log n}{2}\leq \log(n/c_5)\leq \log((1-t)^{-1}) \leq \log(c_5n) \leq 2\log n 
$$ 
since $c_5<2$ and $n\geq 4$. Now 
$$
1-\log(1-t)\sim_3 \log(n)
$$
and hence 
$$
\frac{(1-\log(1-t))^{-\beta}}{(1-t)^2}\sim_{c_2} n^2\log^{-\beta} n, 
$$
where $c_2=3^{\beta} c_5^2$.\\

If $t>e^{-\frac{1}{4}}$, the computations above imply that
$$
\frac{(1-\log(1-t))^{-\beta}}{(1-t)^2}\sim_{c_2} n^2\log^{-\beta} n \sim_{c_1}1+\sum_{k=2}^\infty k\log^{-\beta}(k)t^k,
$$
where $n$ is such that $e^{-\frac{1}{n}}< t \leq e^{-\frac{1}{n+1}}$.
If $0\leq t<e^{-\frac{1}{4}}$, then one has that 
$$
1\leq 1+\sum_{k=2}^\infty k\log^{-\beta}(k)t^k \leq c_\beta,
$$
where the constant $c_\beta>0$ depends only on $\beta$, and 
$$
\frac{1}{d_\beta}\leq \frac{(1-\log(1-t))^{-\beta}}{(1-t)^2}\leq d_\beta, 
$$
where the constant $d_\beta>0$ depends only on $\beta$.
Hence 
$$
\frac{(1-\log(1-t))^{-\beta}}{(1-t)^2}  \sim_c 1+\sum_{k=2}^\infty k\log^{-\beta}(k)t^k  
$$
for all $t\in[0,1)$, where the constant $c\geq 1$ depends on $\beta$.
\end{proof}

\begin{lemma}\emph{\cite[Lemma 3.9]{Geiss}}\label{L3.7}
For $f=\sum_{k=0}^\infty a_kh_k\in L_2(\gamma)$, $t\in [0,1)$ and $Z=f(W_1)$ one has that
\[
\begin{split}
H_WZ^2(t) &= \sum_{k=0}^\infty a_{k+2}^2(k+2)(k+1)t^k,\\
H_SZ^2(t) &= \sum_{k=0}^\infty \left(a_{k+2}-\frac{a_{k+1}}{\sqrt{k+2}}\right)^2(k+2)(k+1)t^k,
\end{split}
\]
where $W$ is a standard Brownian motion and $S$ is the geometric Brownian motion. Moreover 
$$
\frac{1}{12}H_WZ^2(t)-\frac{2}{3}(a_1^2+a_2^2)\leq H_SZ^2(t)\leq 4H_WZ^2(t)+2a_1^2.
$$
\end{lemma}


\section{Application: Approximation of certain d-dimensional stochastic integrals with drift}

We can apply Theorems \ref{main} and \ref{main3} also to the discrete time approximation of $d$--dimensional stochastic integrals 
considered by Zhang \cite{Zhang}, Temam \cite{Temam} and Hujo \cite{Hujo2}. Our setting introduced below recalls for the convenience of the reader line by line the setting of \cite{Hujo2}, which 
generalizes the $1$--dimensional setting of Section 4 to $d$ dimensions.

We assume a stochastic basis $(\Omega,\mathcal{F},\mathbb{P}, (\mathcal{F}_t)_{t\in [0,1]})$, where $(\mathcal{F}_t)_{t\in [0,1]}$ is the augmentation of the natural filtration generated by the $d$--dimensional 
standard Brownian motion $W=(W_t)_{t\in [0,1]}$ with $\mathcal{F}=\mathcal{F}_1$. 

We consider a diffusion $X=(X^1,...,X^d)$, where 
\begin{equation}\label{ddimX}
X_t^i=x_0^i+\int_0^t b_i(X_u)du + \sum_{j=1}^d\int_0^t\sigma_{ij}(X_u)dW_u^j,\quad  t\in[0,1],\  a.s.
\end{equation}
for all $i=1,...,d$ and $x_0=(x_0^1,...,x_0^d)$. We assume $X$ is obtained through $Y$ given 
as unique path--wise continuous solution of 
\begin{equation}\label{sde}
Y_t^i=y_0^i+\int_0^t \hat{b}_i(Y_u)du + \sum_{j=1}^d\int_0^t\hat{\sigma}_{ij}(Y_u)dW_u^j, \quad  t\in[0,1],\  a.s.
\end{equation}
for all $i=1,...,d$, where $\hat{b}_i,\hat{\sigma}_{ij}\in \mathcal{C}_b^\infty(\mathbb{R}^d)$ and 
$(\hat{\sigma}\hat{\sigma}^T)_{ij}(x)=\sum_{k=1}^d\hat{\sigma}_{ik}(x)\hat{\sigma}_{jk}(x)$ 
is uniformly elliptic, i.e. 
$$
\sum_{i,j=1}^d(\hat{\sigma}\hat{\sigma}^T)_{ij}(x)\xi_i\xi_j\geq \lambda\Vert\xi\Vert^2,\ 
\text{ for all }x,\xi\in\mathbb{R}^d \text{ and some }\lambda>0,
$$
where $\Vert\cdot\Vert$ is the Euclidean norm.
Again, we assume that $X$ is obtained through $Y$ by one of the following two ways:
\begin{enumerate}
\item[(a)] $x_0=y_0\in\mathbb{R}^d$, $\hat{b}_i(x):=b_i(x)$, $\hat{\sigma}_{ij}(x):=\sigma_{ij}(x)$, 
and $X_t=Y_t$,
\item[(b)] $x_0=e^{y_0}\in(0,\infty)^d,$ 
$\hat{b}_i(y):=\frac{b_i(e^y)}{e^{y_i}}-\frac{1}{2}\sum_{j=1}^d\hat{\sigma}_{ij}^2(y),$ 
$\hat{\sigma}_{ij}(y):=\frac{\sigma_{ij}(e^y)}{e^{y_i}},$ and $X_t=e^{Y_t}$.
\end{enumerate}
Here and in the following $e^y=(e^{y_1},..., e^{y_d})$ for $y=(y_1,...,y_d)$.
As in one dimensional case, (a) is related to the standard Brownian motion and (b) to 
the geometric Brownian motion. 

Moreover, we assume that $f:E\rightarrow\mathbb{R}$ is a Borel--function such that for some 
$q\in (0, \infty)$ and $C>0$ it holds that
\begin{equation}\label{f<}
|f(x)|\leq C(1+\Vert x\Vert^q),\ x\in E, 
\end{equation}
where the set $E$ is defined by
$$
E:=
\begin{cases}
\mathbb{R}^d\quad&\text{in case }(a)\\
(0,\infty)^d &\text{in case }(b).
\end{cases}
$$
Finally, we define the function $g:\mathbb{R}^d\rightarrow \mathbb{R}$  by
$$
g(y):=
\begin{cases}
f(y)\quad&\text{in case }(a)\\
f(e^y) &\text{in case }(b).
\end{cases}
$$

\begin{theorem}\emph{\cite[Theorem 8 on p. 263]{Friedman1}, \cite[Theorem 5.4 on p. 149]{Friedman2}}\label{app1}
For $\hat{b}$, $\hat{\sigma}$ with $\hat{\sigma}\hat{\sigma}^T$ uniformly elliptic, there exists a transition density 
$\Gamma:(0,1]\times \mathbb{R}^d\times \mathbb{R}^d\rightarrow [0,\infty)\in \mathcal{C}^\infty$ such that 
$\mathbb{P}(Y_t\in B)=\int_B\Gamma(t,y,\xi)d\xi$ for $t\in(0,1]$ and $B\in\mathcal{B}(\mathbb{R}^d)$, where 
$Y=(Y_t)_{t\in[0,1]}$ is the strong solution of stochastic differential equation \eqref{sde} starting in $y$. Moreover, the following is satisfied:
\begin{enumerate}
\item[(i)] For $(s,y,\xi)\in(0,1]\times \mathbb{R}^d\times \mathbb{R}^d$ we have that
$$
\frac{\partial}{\partial s}\Gamma(s,y,\xi)=\frac{1}{2}\sum_{k,l=1}^d\sum_{j=1}^d\hat{\sigma}_{kj}(y)\hat{\sigma}_{lj}(y)
\frac{\partial^2}{\partial y_k\partial y_l}\Gamma(s,y,\xi)+\sum_{i=1}^d \hat{b}_i(y)\frac{\partial}{\partial y_i}\Gamma(s,y,\xi).
$$
\item[(ii)] For $a\in\{0,1,2,...\}$ and multi--indices $b$ and $c$ there exist positive constants $C$ and $D$, depending 
only on $a,b,c$ and $d$ such that
$$
\left|\frac{\partial^{a+|b|+|c|}}{\partial^at\partial^{b}y\partial^{c}\xi}\Gamma(t,y,\xi)\right|\leq 
\frac{C}{t^{(d+2a+|b|+|c|)/2}}e^{-D\frac{\Vert y-\xi\Vert^2}{t}}.
$$ 

\end{enumerate}
  
\end{theorem}

If we apply Theorem \ref{app1} to the stochastic differential equation
$$
\begin{cases}
Z_t^i=Z_0^i+\sum_{j=1}^d \int_0^t\hat{\sigma}_{ij}(Z_u)dW_u^j \quad&\text{in case }(a),\\
Z_t^i=Z_0^i-\int_0^t\left(\frac{1}{2}\sum_{j=1}^d \hat{\sigma}_{ij}^2(Z_u)\right)du+\sum_{j=1}^d \int_0^t\hat{\sigma}_{ij}(Z_u)dW_u^j&\text{in case }(b),
\end{cases}
$$
we obtain a transition density $\Gamma_0$ such that we can define the function 
$G\in \mathcal{C}^\infty([0,1]\times \mathbb{R}^d)$ by 
$$
G(t,y):=\int_{\mathbb{R}^d}\Gamma_0(1-t, y, \xi)g(\xi)d\xi, \ 0\leq t < 1
$$
so that 
$$
\begin{cases}
\left(\frac{\partial}{\partial t} + 
\frac{1}{2}\sum_{k,l=1}^d\left(\hat{\sigma}\hat{\sigma}^T(y)\right)_{kl}\frac{\partial^2}{\partial y_k\partial y_l}\right)G(t,y)=0\ &\text{(a)},\\
\left(\frac{\partial}{\partial t} 
- \sum_{i=1}^d\left(\frac{1}{2}\sum_{j=1}^d \hat{\sigma}_{ij}^2(y)\right)\frac{\partial}{\partial y_i} +
\frac{1}{2}\sum_{k,l=1}^d\left(\hat{\sigma}\hat{\sigma}^T(y)\right)_{kl}\frac{\partial^2}{\partial y_k\partial y_l}\right)G(t,y)=0  &\text{(b)}.
\end{cases}
$$
We define the function $F:E\rightarrow \mathbb{R}$ by setting 
\[
F(t,x):=
\begin{cases}
G(t,x), \quad &\text{ in case }(a),\\
G(t,\log(x)), \quad &\text{ in case }(b),
\end{cases}
\]
where $\log x=(\log(x_1),...,\log(x_d))$, and the operator $\mathcal{L}$ by
$$
\mathcal{L}:=\frac{\partial}{\partial t}+\frac{1}{2}\sum_{k,l=1}^d L_{kl}(x)\frac{\partial^2}{\partial x_k\partial x_l},
$$
where $L_{kl}(x)=\sum_{j=1}^d\sigma_{kj}(x)\sigma_{lj}(x)$. Now we have that 
$$
\mathcal{L}F(t,x)=0 \ \text{ on }[0,1)\times E,
$$
and Itô's formula implies that 
$$
F(t, X_t)=F(0, X_0)+\sum_{k=1}^d\int_0^t\frac{\partial}{\partial x_k}F(u,X_u)dX_u^k,\ \text{a.s. } t\in[0,1).
$$
From Theorem \ref{app1} we get that 
$$
F(t,X_t)\rightarrow f(X_1)\ \text{ in } L_2 \text{ as } t\nearrow 1
$$
and 
$$
f(X_1)=F(0,X_0)+\sum_{k=1}^d\int_0^1 \frac{\partial}{\partial x_k}F(u, X_u)dX_u^k\ \text{a.s.}
$$

\begin{definition}
For $f$, $F$ and $X$ as above we define 
\[
\begin{split}
a_X^{\mathrm{sim}}&(f(X_1), \tau, s)\\&:=\left\Vert\sum_{i=1}^n\sum_{k=1}^d
\int_{t_{i-1}^n\wedge s}^{t_i^n\wedge s}\left(\frac{\partial}{\partial x_k}F(u,X_u)
-\frac{\partial}{\partial x_k}F(t_{i-1}^n,X_{t_{i-1}^n}) \right)dX_u^k\right\Vert_{L_2},
\end{split}
\]
for all $\tau=(t_i)_{i=0}^n\in \mathcal{T}_n$ and $s\in[0,1)$.
\end{definition}

\begin{definition}
We define $H_Xf, H^*_Xf:[0,1)\rightarrow [0,\infty)$ by setting

\[
\begin{split}
&H_Xf(t):=\left(\sup_{\alpha, \beta}\mathbb{E}\left[L_{\alpha\alpha}(X_t)L_{\beta\beta}(X_t)
\left|\frac{\partial^2}{\partial x_\alpha \partial x_\beta}F(t,X_t)\right|^2\right]\right)^\frac{1}{2}\quad \text{and}\\
&H^*_Xf(t):=\sup_{s\in[0,t]}H_Xf(s).
\end{split}
\]

\end{definition}

Finally, we define functions 
$Q_i:\mathbb{R}^d\rightarrow \mathbb{R}$ for $i=1,...,d$ by
$$
Q_i(x):=
\begin{cases}
1,\quad &\text{in case }(a)\\
x_i &\text{in case }(b).
\end{cases}
$$

In this setting we have the following theorem, which refines \cite[Theorem 1]{Hujo2}.

\begin{theorem}\label{a2thm1}
Assume that for all $x\in E$ 
$$
\left|\frac{\partial^s}{\partial_{x_\beta}^q\partial_{x_\alpha}^r}\sigma_{ij}(x)\right|\leq 
C_1\frac{Q_i(x)}{Q_\beta^q(x)Q_\alpha^r(x)},\ \text{ where } q+r=s,\ q,r,s\in\{0,1,2\}, 
$$
$|b_i(x)|\leq C_1Q_i(x)$ and $L_{ii}(x)\geq \frac{1}{C_1}Q_i^2(x)$ for $i\in \{1,...,d\}$ and some fixed $C_1>0$. 
\begin{enumerate}
\item[(1)] If one has that 
\[
I_H:=\int_0^1 H_X^*f(t)dt<\infty,
\]
then 
$$
\inf_{\tau\in \mathcal{T}_n}\sup_{s\in[0,1]}a_X^{\mathrm{sim}}(f(X_1), \tau, s)
\leq\frac{D_1I_H}{\sqrt{n}}\quad\text{for all } n\in\mathbb{N},
$$
where $D_1=D_1(C_1,d)>0$.
\item[(2)] If there exists $C_2>0$ and $\alpha\in(\frac{1}{2},1)$ such that
\[
H_X^*f(t)\leq C_2\frac{(1-\log(1-t))^{-\alpha}}{1-t} \quad \text{for all }t\in[0,1),
\] 
then 
$$
\inf_{\tau\in \mathcal{T}_n}\sup_{s\in[0,1]}a_X^{\mathrm{sim}}(f(X_1), \tau, s)
\leq\frac{D_2}{\sqrt{n^{2\alpha-1}}}\quad\text{for all } n\in\mathbb{N},
$$
where $D_2=D_2(C_1,C_2,\alpha,d)>0$.
\end{enumerate}

\end{theorem}

\begin{proof}[Proof of Theorem \emph{\ref{a2thm1}} ]
Hujo showed in the proof of \cite[Theorem 1, p. 18]{Hujo2} that under the assumptions of Theorem \ref{a2thm1} we have that 
\[
\begin{split}
&\mathbb{E}\left|\sum_{i=1}^n\sum_{k=1}^d
\int_{t_{i-1}^n\wedge s}^{t_i^n\wedge s}\left(\frac{\partial}{\partial x_k}F(u,X_u)
-\frac{\partial}{\partial x_k}F(t_{i-1}^n,X_{t_{i-1}^n}) \right)dX_u^k\right|^2\\
&\leq c\sum_{i=1}^n\int_{t_{i-1}^n}^{t_i^n}\int_{t_{i-1}^n}^t \sup_{\alpha, \beta}\mathbb{E}\left[L_{\alpha\alpha}(X_u)L_{\beta\beta}(X_u)
\left|\frac{\partial^2}{\partial x_\alpha \partial x_\beta}F(u,X_u)\right|^2\right]dudt\\
&\leq c\sum_{i=1}^n\int_{t_{i-1}^n}^{t_i^n}(t_i-t)[H_X^*f(t)]^2dt
\end{split}
\]
for any  $s\in[0,1)$ and any time net $\tau=(t_i^n)_{i=0}^n$, where $c=c(C_1, d)$. Hence we can conclude by Theorems \ref{main} and \ref{main3}. 
\end{proof}


\subsection*{Acknowledgements} The author wishes to express his gratitude to his supervisor Stefan Geiss for suggestions concerning 
the contents of this paper, and for reading the manuscript. The author would also like to thank Eero Saksman for his support as well as Academy of Finland project no. 7110599, Finnish Graduate School of Stochastics and Statistics, the Finnish Academy of Science and Letters, Vilho, Yrj\"o, and Kalle V\"ais\"al\"a Foundation, and Magnus Ehrnrooth Foundation for their financial support.



\end{document}